\documentstyle{article}\begin{document}\centerline{\bf A remarkable Definite Integral}\vskip .2in

\centerline{ M.L. Glasser}\vskip .2in
\centerline{Physics Department, Clarkson University}
\centerline{ Potsdam, NY 13699-5820 }\vskip 1in
\centerline{\bf Abstract}\vskip .1in
\begin{quote}
A definite integral is evaluated explicitly whose integrand, subject to mild restrictions, contains an arbitrary function. A number of example are given along with a new integral of the Riemann zeta function.

\end{quote}

\newpage

\section{Introduction}
The ultimate goal in  the definite integration  of functions of a single variable would be the formula
$$\int_{-\infty}^{\infty}F(x)dx=G[F]\eqno(0)$$
where where $F$ denotes any integrable function and  $G[F]$ is an explicit expression free of the integration symbol. A trivial, but more restricted,  example is the definition of the Dirac delta function $\int f(x)\delta(x-a)dx=f(a)$.  This note deals with a less trivial, but, to the author, quite remarkable formula of this nature and some of its consequences.
\vskip .2in
\section{Calculation}

We start with the formula proven in [1]: For $a>0$, $t>0$

$$\int_0^{\infty}\frac{e^{-t x^2}\cos(t\pi x)\cosh x}{1+2a^2\cosh(2x)+a^4}dx=\frac{\pi e^{-t[\frac{\pi^2}{4}+(\ln a)^2]}}{4a(1+a^2)}.\eqno(1)$$
Now, multiply both sides of (1) by any function $f(t)$, which, for simplicity, we shall take to be real for real argument, possessing a Laplace transform
$$F(k)=\int_0^{\infty}e^{-kt}f(t) dt\eqno(2)$$
and integrate with respect to $t$ over $[0,\infty]$.  If, in addition to (2), we require that interchange of the order of integration be allowed,  one has the identity
$$\int_{-\infty}^{\infty}\frac{F(x^2+i\pi x)\cosh x}{1+2a^2\cosh(2x)+a^4} dx
=\frac{\pi F[(\pi/2)^2+(\ln a)^2]}{2a(1+a^2)}.\eqno(3)$$
Equation (3) is our principal result. This appears to be the first formula of this nature (although one might consider Ramanujan's  ``Master Theorems" [2]  as belonging to this class) and, in view of the mild restrictions placed of the function $F$, (3) is nearly unlimited in  scope. 
\vskip .2in

\section{ Examples}

For $f(t)=e^{-bt}$, one has $F(k)=1/(b+k)$, $b>0$. Hence
$$I(a,b)=\int_0^{\infty}\frac{(x^2+b)\cosh x}{(x^4+(2b+\pi^2)x^2+b^2)(1+2a^2\cosh(2x)+a^4)}dx$$
$$=\frac{\pi}{4a(1+a^2)(b+\pi^2/4+\ln^2 a)}.\eqno(4)$$
E.g. for $a=0.7$, $b=2$, (4) gives 0.163891=0.163891. On the other hand, $\lim_{b\to 0}I(a,b)\ne I(a,0)$.\vskip .1in

For $f(t)=J_0(t)$, $F(k)=1/\sqrt{1+k^2}$. Hence
$$\int_{-\infty}^{\infty}\frac{\cosh x}{1+2a^2\cosh(2x)+a^4}\frac{dx}{\sqrt{1+(x^2+\pi i x)^2}}$$
$$=\frac{\pi}{2a(1+a^2)\sqrt{1+(\frac{\pi^2}{4}+\ln^2a)^2}}.\eqno(5)$$
E.g. for $a=0.7$, (5) gives  0.000708622 on both sides.\vskip .3in

For $F(k)=e^{-bk^2}$ one obtains
$$\int_0^{\infty}e^{-bx^2(x^2-\pi^2)}\frac{\cos(2b\pi x^3)\cosh x}{1+2a^2\cosh(2x)+a^4}dx=e^{-b(\pi^2/4+\ln^2a)^2}\frac{\pi}{4a(1+a^2)}.\eqno(6)$$ For example, with $a=b=0.3$ both sides of (6) give 0.0240764.

 Another case in which the real part can be made explicit is $F(k)=\cos(\alpha k)$, which leads to
 $$\int_0^{\infty}\frac{\cos(\alpha x^2)\cosh(\alpha\pi x)\cosh x}{1+2a^2\cosh(2x)+a^4}dx=\frac{\pi \cos[\alpha(\pi^2/4+\ln^2a)]}{4a(1+a^2)}\eqno(7)$$
 as can be verified numerically for $\alpha\pi\le1$. E.g. for $\alpha=0.1$ and $a=1+2i$, both sides of (7) give $-0.0783703+0.00264214i$. It is interesting that in this case $F(k)$ is the Laplace transform of a distribution.\vskip .2in
 
 \section{Discussion}
 Although at the present time (3) appears to be unique, all that is required is an integral evaluation where the same integral transform kernel appears under the integral sign on the left and "naked" on the right. For example, (6) and (7) are obvious candidates where the former involves the Laplace transform(wrt b) and the later the cosine transform (wrt $\alpha$). Neither of these formulas produces anything new.
 
 By setting $a=1$ in (3) one has
$$\int_{-\infty}^{\infty}F(x^2+i\pi x){\rm sech} x\; dx=\pi F(\pi^2/4).\eqno(8)$$
For example,
$$\int_{\infty}^{\infty}\frac{dx}{\cosh (\pi x) \Gamma(4ax(x+i )+b)}=\frac{1}{\Gamma(a+b)}\eqno(9)$$
for real $a,b$. 

One also has, for $a>0$ and
$$\sigma_{2l}=\sum_{n_1<n_2<\cdots<n_{2l}}\prod_{k\ne n_1,\dots,n_{2l}}(x^2+a+\frac{k}{\pi^2})$$
$$\int_{-\infty}^{\infty}\frac{\sum_{l=0}^{\infty}(-1)^l\sigma_{2l} x^{2l}+\left\{\begin{array}{c}
(-1)^{n/2}x^n; n\; even\\
0; n\; odd
\end{array}\right\}}{\prod_{k=0}^{n-1}[(x^2+a+\frac{k}{\pi^2})^2+x^2]}\frac{dx}{\cosh(\pi x)}$$
$$=\prod_{k=0}^{n-1}\left(a+\frac{1}{4}+\frac{k}{\pi^2}\right)^{-1}.\eqno(10)$$

A particularly important result relating to the Riemann Zeta function is
$$\int_{-i\infty}^{i\infty}\frac{x^{s(1-s)}}{\cos(\pi s)\zeta^n[4a s(1-s)]}\frac{ds}{2\pi i}=\frac{x^{1/4}}{2\pi\zeta^n(a)}\eqno(12)$$
which is valid for all  $a>0$, $|x|<1$ and $n=0,1,2, 3,4.$

\vskip .2in
\noindent
{\bf Acknowlegement}
The author thanks Drs. Michael Milgram and Jan Grzesik for their comments on this material.

.\vskip .2in
\centerline{\bf Reference}\vskip .1in

\noindent
[1] M.L. Glasser, {\it Generalization of a Definite Integral of Ramanujan, } J. Ind. Math. Soc.
37, 351 (1974).

\noindent
[2] G.H. Hardy, {\it Ramanujan},[Chelsea Publishers, NY(1940); p.186]

\end{document}